\documentclass[12pt]{article}

\usepackage{amssymb}
\usepackage{amsmath,amsthm}
\usepackage[dvips]{graphicx}
\usepackage{wrapfig}
\usepackage{url}

\newtheorem{theorem}{Theorem}[section]

\newtheorem{lemma}[theorem]{Lemma}
\newtheorem{proposition}[theorem]{Proposition}

\theoremstyle{definition}

\newtheorem{example}[theorem]{Example}
\theoremstyle{remark}

\title{The warping matrix of a knot diagram}
\author{ Ayaka Shimizu \thanks{Department of Mathematics, Gunma National College of Technology, 580 Toriba-cho, Maebashi-shi, Gunma, 371-8530 Japan. shimizu@nat.gunma-ct.ac.jp}}
\date{\today}

\begin{document}

\maketitle

\begin{abstract}
We introduce the warping matrix which is a new description of oriented knots from a viewpoint of warping degree. 
\end{abstract}

\section{Introduction}

Warping degree defined by Kawauchi in \cite{kawauchi} represents a complexity of an oriented knot diagram, and has been studied for knots, links and spatial-graphs. 
\begin{figure}[h]
\begin{center}
\includegraphics[width=130mm]{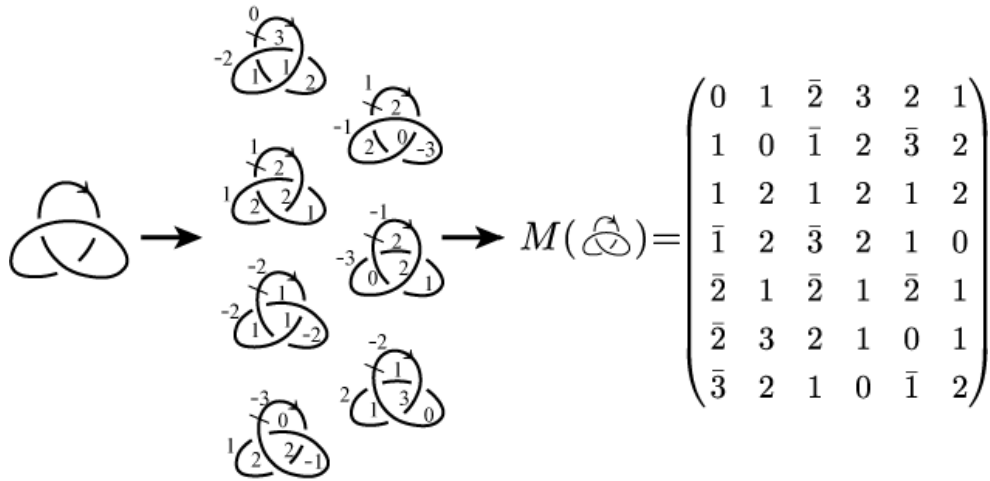}
\caption{Warping matrix of an oriented knot diagram.}
\label{image}
\end{center}
\end{figure}
In this paper, we define the {\it warping matrix} $M(D)$ (resp. $M(P)$) of an oriented knot diagram $D$ (resp. projection $P$) as depicted in Fig. \ref{image}, and show the following theorem: 

\phantom{x}
\begin{theorem}
The warping matrix of an oriented knot diagram represents the oriented knot diagram uniquely. 
\label{mainthm}
\end{theorem}
\phantom{x}

\noindent The rest of this report is organized as follows: 
In Section 2, we define the warping matrix of an oriented knot projection, and look into properties. 
In Section 3, we define the warping matrix of an oriented knot diagram, and prove Theorem \ref{mainthm}. 
In Appendix, we consider a puzzle as an application.

\section{Warping matrix for a knot projection}
\label{w-matrix-projection}

\begin{figure}[h]
\begin{center}
\includegraphics[width=30mm]{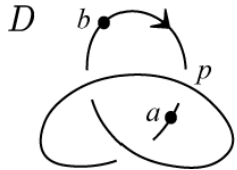}
\caption{The crossing $p$ is a warping crossing point of $D_b$. }
\label{wcp-ex}
\end{center}
\end{figure}

In this section we review the warping degree, and define the warping matrix of an oriented knot projection. 
Let $D$ be an oriented knot diagram on $S^2$. 
Take a base point $b$ of $D$. 
We denote by $D_b$ the based diagram.  
Go along $D$ with the orientation from $b$ to $b$. 
Then we encounter every crossing twice -- as an overcrossing  once and undercrossing once. 
We say that a crossing is a {\it warping crossing point of} $D_b$ if we meet the crossing as an undercrossing first (and overcrossing later). 
For example, the crossing $p$ of the diagram $D$ in Fig. \ref{wcp-ex} is a warping crossing point of $D_b$, and is not a warping crossing point of $D_a$. 
The {\it warping degree $d(D_b)$ of $D_b$} is the number of warping crossing points of $D_b$. 
For example, we have $d(D_b)=2$ and $d(D_a)=1$ in Fig. \ref{wcp-ex}. 
The following lemma was shown in \cite{shimizu}: 

\phantom{x}
\begin{lemma}{$($Lemma 2.5 in \cite{shimizu}$)$}
For two base points $a_1,a_2$ $($resp. $b_1,b_2)$ which are placed on opposite side of an overcrossing  $($resp. undercrossing$)$ as shown in Figure \ref{lem25}, 
we have $d(D_{a_2})=d(D_{a_1})+1$ $($resp. $d(D_{b_2})=d(D_{b_1})-1)$. 
\begin{figure}[h]
\begin{center}
\includegraphics[width=70mm]{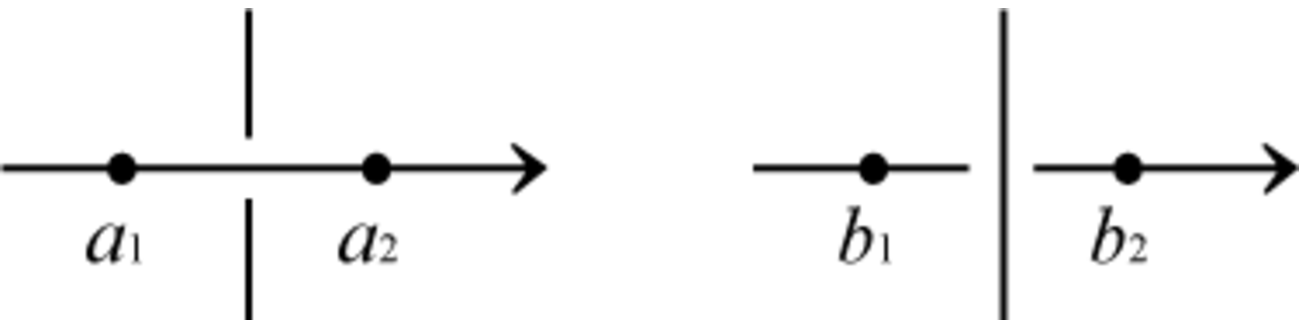}
\caption{Warping degree increases (or decreases) by one when we pass a crossing.}
\label{lem25}
\end{center}
\end{figure}
\label{lemma25}
\end{lemma}
\phantom{x}

\begin{figure}[h]
\begin{center}
\includegraphics[width=30mm]{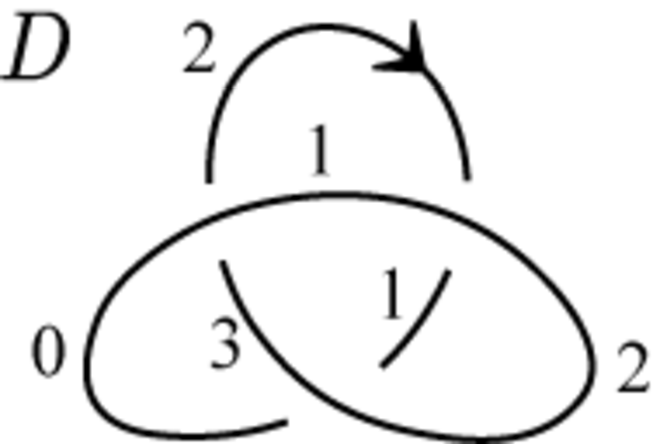}
\caption{Warping degree labeling. }
\label{labeling-ex}
\end{center}
\end{figure}

\noindent An {\it edge} of a knot diagram is a path between crossing points which has no crossings in the interior. 
{\it Warping degree labeling} to $D$ is the labeling giving $d(D_b)$ to each edge $e$, where $b$ is a base point on $e$ (\cite{shimizu-poly}). 
An example is shown in Fig. \ref{labeling-ex}. 
Go along $D$ and read all the warping degree labels from a label. 
Thus we obtain a {\it warping degree sequence}. 
For example, 210123, 101232 and 012321 are warping degree sequences of $D$ in Fig. \ref{labeling-ex}. 

Now we define the warping matrix. 
Let $P$ be an oriented knot projection on $S^2$ with $c$ crossings. 
We obtain $2^c$ knot diagrams from $P$ by giving over/under information to each crossing. 
Consider a $2^c \times 2c$ matrix such that each row represents the warping degree sequence starting from the same edge for all the $2^c$ knot diagrams. 
We call such a matrix a {\it warping matrix $M(P)$ of $P$}. 
For example, we have the following warping matrices:

\newsavebox{\boxa}
\newsavebox{\boxb}

\sbox{\boxa}{\includegraphics[width=11mm]{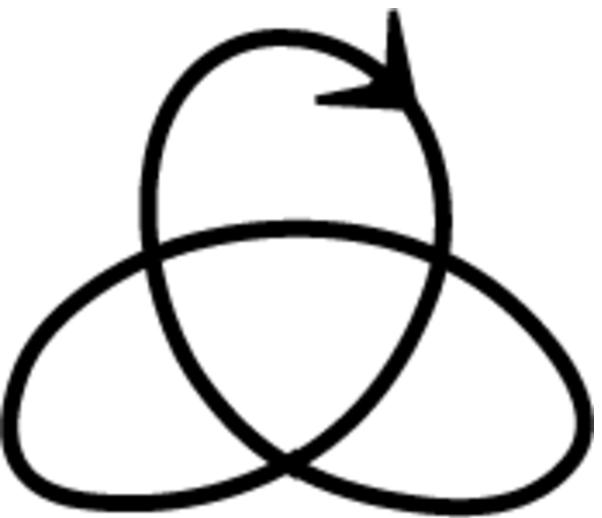}}
\sbox{\boxb}{\includegraphics[width=13mm]{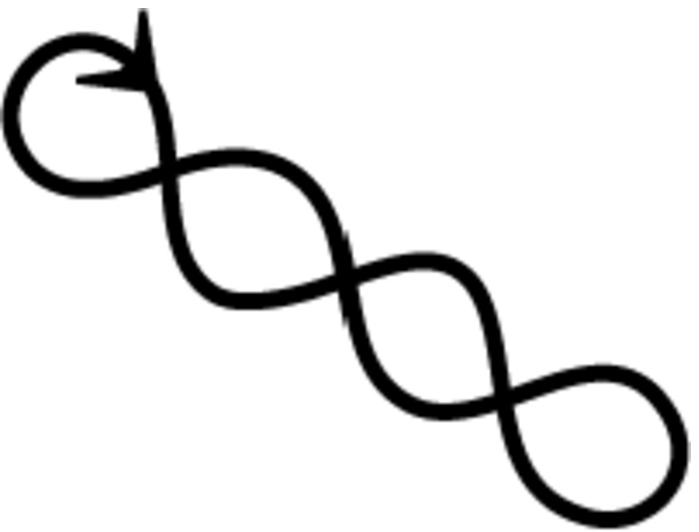}}

\begin{equation*}
M( \parbox[c]{11mm}{\usebox{\boxa}})=
\begin{pmatrix}
0 & 1 & 2 & 3 & 2 & 1 \\
1 & 0 & 1 & 2 & 3 & 2 \\
1 & 2 & 1 & 2 & 1 & 2 \\
1 & 2 & 3 & 2 & 1 & 0 \\
2 & 1 & 0 & 1 & 2 & 3 \\
2 & 1 & 2 & 1 & 2 & 1 \\
2 & 3 & 2 & 1 & 0 & 1 \\
3 & 2 & 1 & 0 & 1 & 2 \\
\end{pmatrix}
, \ M( \parbox[c]{13mm}{\usebox{\boxb}})=
\begin{pmatrix}
0 & 1 & 2 & 3 & 2 & 1 \\
1 & 0 & 1 & 2 & 1 & 0 \\
1 & 2 & 1 & 2 & 1 & 2 \\
1 & 2 & 3 & 2 & 3 & 2 \\
2 & 1 & 0 & 1 & 0 & 1 \\
2 & 3 & 2 & 1 & 2 & 3 \\
2 & 1 & 2 & 1 & 2 & 1 \\
3 & 2 & 1 & 0 & 1 & 2 \\
\end{pmatrix}
\end{equation*}

\noindent From a knot projection, we have warping matrices not uniquely, and they are equivalent up to the following two moves:

\begin{itemize}
\item[(R)] Switch two rows. 
\item[(C)] Apply a cyclic permutation on columns. 
\end{itemize}

\noindent (R) means we have $2^c!$ choices of the order of knot diagrams and (C) means we have $2c$ choices of the start point. 
We consider warping matrices up to those moves. 
We have the following proposition:

\phantom{x}

\begin{proposition}
Let $P$ be an oriented knot projection with $c$ crossings. 
A warping matrix $M(P)=(a_{i j})$ satisfies the following (i)--(v): 

\begin{itemize}
\item[(i)] $\vert a_{i (j+1)}-a_{i j} \vert = \vert a_{i 1}-a_{i (2c)} \vert =1$ for any $i$ and $j$. 
\item[(ii)] At each column, the number $n$ appears 
$\begin{pmatrix}
c \\
n
\end{pmatrix}$
times $(n=0,1,2, \dots , c)$. 
\item[(iii)] There are just $2^{c-1}$ disjoint pairs of rows uniquely such that the sum of them is $(c \ c \dots c)$. 
\item[(iv)] There are just two rows $( k \ k+1 \ k \ k+1 \ k \dots k+1 )$ and $(l \ l-1 \ l \ l-1 \ l \dots l-1 )$, where $k+l=c$. 
\end{itemize}

\begin{proof}
(i): By the property of warping degree labeling (see Lemma \ref{lemma25}). 
(ii): See the proof of Theorem 1.1 in \cite{KS}. 
(iii): Let $D$ be an oriented knot diagram, and $D^*$ the diagram obtained from $D$ by applying crossing changes at all the crossings of $D$. 
We have $d(D_b)+d({D^*}_b)=c(D)$ for any base point $b$ (see Lemma 2.1 in \cite{shimizu} and Example 2.4 in \cite{shimizu}), and we have $2^{c-1}$ pairs of knot diagrams such as $D$ and $D^*$ from $P$. 
(iv): We have just two alternating diagrams from $P$. 
\end{proof} 
\label{rule}
\end{proposition}

\phantom{x}

\noindent In the following example, we consider all $2 \times 2$, $4 \times 4$ and $8 \times 6$ matrices satisfying the properties of Proposition \ref{rule}. 

\phantom{x}
\begin{example}
Any $2^c \times 2c$ matrix satisfying (i)--(iv) in Proposition \ref{rule} is equivalent up to the moves (R) and (C) and the vertical reflection to 
\begin{equation*}
\begin{pmatrix}
0 & 1 \\
1 & 0 \\
\end{pmatrix}
 \ {\rm for} \ c=1, 
 \end{equation*}
 
\begin{equation*}
\begin{pmatrix}
0 & 1 & 2 & 1 \\
1 & 0 & 1 & 0 \\
1 & 2 & 1 & 2 \\
2 & 1 & 0 & 1 \\
\end{pmatrix}
 \ {\rm for} \ c=2, 
\end{equation*}
and 
\begin{equation*}
\begin{pmatrix}
0 & 1 & 2 & 3 & 2 & 1 \\
1 & 0 & 1 & 2 & 3 & 2 \\
1 & 2 & 1 & 2 & 1 & 2 \\
1 & 2 & 3 & 2 & 1 & 0 \\
2 & 1 & 0 & 1 & 2 & 3 \\
2 & 1 & 2 & 1 & 2 & 1 \\
2 & 3 & 2 & 1 & 0 & 1 \\
3 & 2 & 1 & 0 & 1 & 2 \\
\end{pmatrix}
, \ 
\begin{pmatrix}
0 & 1 & 2 & 3 & 2 & 1 \\
1 & 0 & 1 & 2 & 1 & 0 \\
1 & 2 & 1 & 2 & 1 & 2 \\
1 & 2 & 3 & 2 & 3 & 2 \\
2 & 1 & 0 & 1 & 0 & 1 \\
2 & 3 & 2 & 1 & 2 & 3 \\
2 & 1 & 2 & 1 & 2 & 1 \\
3 & 2 & 1 & 0 & 1 & 2 \\
\end{pmatrix}
 {\rm or}  
\begin{pmatrix}
0 & 1 & 2 & 1 & 2 & 1 \\
1 & 0 & 1 & 0 & 1 & 0 \\
1 & 2 & 1 & 2 & 3 & 2 \\
1 & 2 & 3 & 2 & 1 & 2 \\
2 & 1 & 0 & 1 & 2 & 1 \\
2 & 3 & 2 & 3 & 2 & 3 \\
2 & 1 & 2 & 1 & 0 & 1 \\
3 & 2 & 1 & 2 & 1 & 2 \\
\end{pmatrix}
\end{equation*}
for $c=3$. 
\end{example}

\phantom{x}

\begin{figure}[ht]
\begin{center}
\includegraphics[width=100mm]{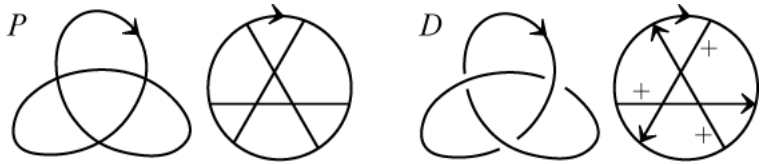}
\caption{Gauss diagrams. }
\label{gauss}
\end{center}
\end{figure}

Next, we review Gauss diagrams. 
Kauffman introduced Gauss code in \cite{kauffman}, and Goussarov, Polyak and Viro represented Gauss codes visually by Gauss diagrams in \cite{GPV}. 
Let $P$ be an oriented knot projection on $S^2$. 
Now we consider $P$ as an immersion $P: S^1 \to S^2$ of the circle into the sphere with some double points (crossings). 
A {\it Gauss diagram} for $P$ is an oriented circle considered as the preimage of the immersed circle with chords connecting the preimages of each crossing. 
Let $D$ be an oriented knot diagram on $S^2$. 
We obtain the Gauss diagram for $D$ in the same way as knot projections by giving the orientation from overcrossing to undercrossing and the crossing sign to each chord (see Fig. \ref{gauss}). 
A Gauss diagram for a knot diagram represents the knot diagram uniquely. 
We have the following lemma: 

\phantom{x}
\begin{lemma}
A warping matrix of an oriented knot projection represents a Gauss diagram for an oriented knot projection uniquely. 

\begin{proof}
Let $M$ be a $2^c \times 2c$ matrix which is a warping matrix of a knot projection, where $c$ is a positive integer. 
Let $A$ be the $2c \times 2c$ matrix defined by 
\begin{equation*}
A= \left(
\begin{array}{ccccc}
-1 & 0 & \ldots & 0 & 1 \\
1 & -1 & \ldots & 0 & 0 \\
0 & 1 & \ldots & 0 & 0 \\
\vdots & \vdots & \ddots & \vdots & \vdots \\
0 & 0 & \ldots & -1 & 0 \\
0 & 0 & \ldots & 1 & -1 \\
\end{array}
\right) 
\end{equation*}

\noindent and let $U=MA$. 
Each element of $U$ is 1 or -1 because of Lemma \ref{lemma25}, and $U$ is a matrix such that each row represents an $ou$ sequence (see \cite{nakanishi, shimizu-link}) where 1 implies $o$ and -1 implies $u$. 
Hence each column of $U$ corresponds to a crossing. 
We can divide the $2c$ columns of $U$ into $c$ pairs such that the sum of the two columns is ${\bf 0}$ because for each column, there exists a column such that their sum is ${\bf 0}$ (they represents the same crossing) 
and there do not exist the same two columns (a warping matrix has all over/under information). 
Thus we have the correspondence of columns representing the same crossing uniquely, and obtain the Gauss diagram for a knot projection. 
\end{proof}
\end{lemma}

\noindent Here is an example. 

\begin{example}

For
\begin{equation*}
M= \left(
\begin{array}{cccccc}
0 & 1 & 2 & 1 & 2 & 1 \\
1 & 0 & 1 & 0 & 1 & 0 \\
1 & 2 & 1 & 2 & 3 & 2 \\
1 & 2 & 3 & 2 & 1 & 2 \\
2 & 1 & 0 & 1 & 2 & 1 \\
2 & 3 & 2 & 3 & 2 & 3 \\
2 & 1 & 2 & 1 & 0 & 1 \\
3 & 2 & 1 & 2 & 1 & 2 \\
\end{array}
\right) 
\end{equation*}

\noindent we have
\begin{equation*}
U=MA= \left(
\begin{array}{cccccc}
0 & 1 & 2 & 1 & 2 & 1 \\
1 & 0 & 1 & 0 & 1 & 0 \\
1 & 2 & 1 & 2 & 3 & 2 \\
1 & 2 & 3 & 2 & 1 & 2 \\
2 & 1 & 0 & 1 & 2 & 1 \\
2 & 3 & 2 & 3 & 2 & 3 \\
2 & 1 & 2 & 1 & 0 & 1 \\
3 & 2 & 1 & 2 & 1 & 2 \\
\end{array}
\right) 
\left(
\begin{array}{cccccc}
-1 & 0 & 0 & 0 & 0 & 1 \\
1 & -1 & 0 & 0 & 0 & 0 \\
0 & 1 & -1 & 0 & 0 & 0 \\
0 & 0 & 1 & -1 & 0 & 0 \\
0 & 0 & 0 & 1 & -1 & 0 \\
0 & 0 & 0 & 0 & 1 & -1 \\
\end{array}
\right) 
\end{equation*}
\begin{equation*}
=\left(
\begin{array}{cccccc}
1 & 1 & -1 & 1 & -1 & -1 \\
-1 & 1 & -1 & 1 & -1 & 1 \\
1 & -1 & 1 & 1 & -1 & -1 \\
1 & 1 & -1 & -1 & 1 & -1 \\
-1 & -1 & 1 & 1 & -1 & 1 \\
1 & -1 & 1 & -1 & 1 & -1 \\
-1 & 1 & -1 & -1 & 1 & 1 \\
-1 & -1 & 1 & -1 & 1 & 1 \\
\end{array}
\right)
\end{equation*}

\noindent and obtain the three pairs of columns 1st and 6th, 2nd and 3rd, and 4th and 5th columns whose sum is ${\bf 0}$. 
Thus we obtain the Gauss diagram in Fig. \ref{gauss-ex}. 
\begin{figure}[h]
\begin{center}
\includegraphics[width=80mm]{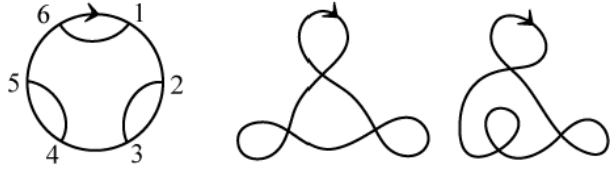}
\caption{We obtain a Gauss diagram from $M$. }
\label{gauss-ex}
\end{center}
\end{figure}

\label{ex-wm-p}
\end{example}

\section{Proof of Theorem \ref{mainthm}}
\label{w-matrix-diagram}

\begin{figure}[h]
\begin{center}
\includegraphics[width=100mm]{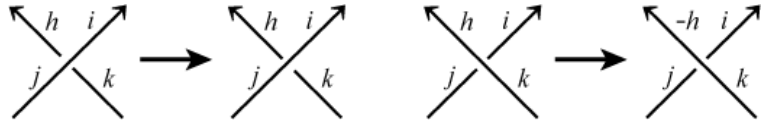}
\caption{Add a minus. }
\label{add-signs}
\end{center}
\end{figure}

\begin{figure}[h]
\begin{center}
\includegraphics[width=30mm]{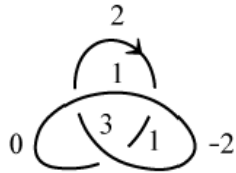}
\caption{Signed warping degree labeling. }
\label{signed-label}
\end{center}
\end{figure}

In this section, we define the warping matrix for oriented knot diagrams. 
First, we define a signed warping degree sequence. 
For an oriented knot diagram $D$ with warping degree labeling, we add signs as follows: 
Go along $D$, and if we encounter a negative crossing as overcrossing, add a minus to the label after the crossing as depicted in Fig. \ref{add-signs}. 
Thus we obtain the {\it signed warping degree labeling}, and obtain a {\it signed warping degree sequence} from it. 
For example, $2101\bar{2}3$ is a signed warping degree sequence of the diagram in Fig. \ref{signed-label}. 

Let $D$ be an oriented knot diagram on $S^2$ with $c$ crossings, and $P$ the knot projection obtained from $D$ by forgetting the over/under information. 
Let $M(P)$ be a warping matrix of $P$, and $m(P)$ the matrix obtained from $M(P)$ by replacing each row with {\it signed} warping degree sequence. 
Note that for each column, there are no or  $2^{c-1}$ bars. 
Let $M(D)$ be the $(2^c -1) \times 2c$ matrix obtained from $m(P)$ by deleting the row representing a signed warping degree sequence of $D$. 
We call $M(D)$ the {\it warping matrix of $D$}, and consider it up to the moves (R) and (C). 
We prove Theorem \ref{mainthm}: 

\phantom{x}
\noindent {\it Proof of Theorem \ref{mainthm}} \ 
From $M(D)$, we can restore $M(P)$ by Proposition \ref{rule} (ii), 
and we can restore the signed warping degree sequence of $D$ by counting minus at each column of $M(D)$. 
Let $G$ be the Gauss diagram of $M(P)$. 
We can give orientations and signs to the chords of $G$ by the signed warping degree sequence of $D$. 
Thus we obtain a Gauss diagram, and it represents an oriented knot diagram uniquely. 
\hfill$\square$
\phantom{x}

\section*{Appendix}

In this appendix, we introduce a Sudoku-like puzzle as an application of the warping matrix of a knot projection. 
At first, there is a $2^c \times 2c$ grid filled with some digits initially, where $c$ is a positive integer. 
The objective of this puzzle is to fill the grid so that the placement of digits satisfies the rule of Proposition \ref{rule} (or just (i) and (ii) in Proposition \ref{rule} for simplicity).  
Here is an example:

\begin{center}
\begin{tabular}{c}

\begin{minipage}{60mm}
\begin{center}
\begin{tabular}{|c|c|c|c|c|c|}\hline
  &  &  &  & 3 & \\ \hline
  &  & 3 &  &  & \\ \hline
  &  & 2 &  & 0 & \\ \hline
 3 &  &  & 0 &  & \\ \hline
  & 1 &  & 1 &  & 1 \\ \hline
  & 2 &  & 2 &  & 2 \\ \hline
 0 &  &  & 3 &  & \\ \hline
  &  & 0 &  &  & \\ \hline
\end{tabular}
\end{center}
\end{minipage}

\begin{minipage}{80mm}
\begin{center}
\begin{tabular}{|c|c|c|c|c|c|c|c|}\hline
 2 & & & & & 3 & &  \\ \hline  
 & 4 & & & & 0 & &  \\ \hline 
1 & & 3 & & 1 & & &  \\ \hline  
 & & & 0 & & & &  \\ \hline 
4 & & & & 2 & & &  \\ \hline 
1 & & & & & 4 & &  \\ \hline 
 & & & 1 & & & 4 &  \\ \hline 
 & & & & & & 1 &  \\ \hline  
 & & 0 & & & 3 & & 3  \\ \hline  
1 & 2 & & & 1 & & 1 &  \\ \hline   
 & & & 3 & 2 & 3 & &  \\ \hline 
 & & 3 & 2 & & & & 2  \\ \hline 
 & 3 & & 3 & & & &  \\ \hline  
 & 3 & & & & & & 3 \\ \hline 
 & & 3 & & & & & 0 \\ \hline 
 & & 4 & & & & 2 &  \\ \hline 
\end{tabular}
\end{center}
\end{minipage}

\end{tabular}
\end{center}

\noindent The left grid represents a standard projection of a trefoil knot, and the right one represents the standard projection of a figure-eight knot.

\bibliographystyle{amsalpha}

\end{document}